\documentclass[12pt]{amsart}

\usepackage{amsfonts}
\usepackage{amsmath}
\usepackage{amssymb}
\usepackage{amscd}
\usepackage[utf8]{inputenc}
\usepackage{graphicx}
\usepackage{enumitem}
\usepackage{mdwlist}
\usepackage[noend]{algorithmic}

\theoremstyle{plain}
\newtheorem{theorem}{Theorem}[section]

\theoremstyle{definition}

\newtheorem{algo}[theorem]{Algorithm}

\theoremstyle{remark}
\newtheorem{remark}[theorem]{Remark}
\newtheorem{example}[theorem]{Example}

\begin{document}

\title[Second-Species Counterpoint]{A Projection-Oriented Mathematical Model for Second-Species Counterpoint}
\author{Octavio A. Agust\'{\i}n-Aquino \and Guerino Mazzola}
\address{Instituto de Física y Matemáticas, Universidad Tecnológica de la Mixteca, Huajuapan de León, Oaxaca, México}
\email{octavioalberto@mixteco.utm.mx}
\address{School of Music, University of Minnesota, MN, USA}
\email{mazzola@umn.edu}
\date{September 28th, 2018}
\thanks{This work was partially supported by a grant from the \emph{Niels Hendrik Abel Board}.}

\begin{abstract}
Drawing inspiration from both the classical Guerino Mazzola's symmetry-based model for first-species
counterpoint (one note against one note) and Johann Joseph Fux's \emph{Gradus ad Parnassum}, we propose an extension for second-species (two notes against one note).
\end{abstract}

\keywords{Second-species, counterpoint}
\subjclass[2010]{00A65,05E18}

\maketitle

\section{Introduction}

Guerino Mazzola's counterpoint model, founded on the concepts of
\begin{enumerate}
\item \emph{strong dichotomy}, which encodes the notion of consonance and dissonance, and
\item \emph{counterpoint symmetry}, which is the carrier of
contrapuntal tension and allows to deduce the rules of counterpoint,
\end{enumerate}
has been successful in explaining the necessity of regarding the fourth as a dissonance and obtaining
the general prohibition of parallel fifths as a theorem. It also allows to define new understandings of consonance and dissonance,
thereby leading to the concept of \emph{counterpoint world}, i.e., paradigms for the handling of two-voice
compositions represented as digraphs, whose vertices are consonant intervals and an arrow connects two of them
whenever we have a valid progression. This, in turn, allows us to \emph{morph} one world into another.
See the monograph \cite{AJM15} and the treatise \cite[Part VII]{gM17} for a thorough account.

Despite these accomplishments, Mazzola's model is restricted to the case of \emph{first-species} counterpoint,
which means that only one note can be placed against another. Hence, in order to increase the potential of Mazzola's
model for analysis and composition, it is indispensable to extend it to \emph{second-species} counterpoint (i.e.,
two notes against one) and further. Our approach for a first step in this direction is to extend the notion
of counterpoint interval to a $2$-interval, i.e., one such that two intervals are attached to a cantus firmus, the first
one coming in the downbeat and the second one in the upbeat.

For our extension, the main idea is that the counterpoint symmetries in this case do not determine
another $2$-interval successor, but a first-species interval in the downbeat. The idea behind
this is to blend the species of counterpoint more easily.

\section{General Overview of Mazzola's Counterpoint Model}

Here we quickly survey the key aspects of Mazzola's counterpoint model (we refer the 
reader to \cite{AJM15} and \cite[Part VII]{gM17} for a complete account). We consider the action of the group
\[
 \overrightarrow{GL}(\mathbb{Z}_{2k}) := \mathbb{Z}_{2k}\rtimes \mathbb{Z}_{2k}^{\times}
\]
(which we call the group of \emph{general affine symmetries}) on $\mathbb{Z}_{2k}$,
which can be described in the following manner:
\[
 T^{u}.v(x) = vx+u;
\]
here $T^{u}$ is the \emph{transposition} by $u$, and $v$ is the \emph{linear part} of
the transformation.

We know \cite{oA09,AJM15} that, for any
$k>4$, there is at least one dichotomy $\Delta=(X/Y)$ of $\mathbb{Z}_{2k}$ such that
there is a unique $p\in \overrightarrow{GL}(\mathbb{Z}_{2k})$ and
\[
 p(X)=Y\quad\text{and}\quad p\circ p = \mathrm{id}_{\mathbb{Z}_{2k}},
\]
which is called the \emph{polarity} of the dichotomy. The dichotomies with
this property are called \emph{strong}, and represent the division of intervals
into generalized \emph{consonances} $X$ and \emph{dissonances} $Y$.

Next we consider the \emph{dual numbers}
\[
 \mathbb{Z}_{2k}[\epsilon] = 
 \frac{\mathbb{Z}_{2k}[\mathcal{X}]}{\langle \mathcal{X}^{2}\rangle}=
 \{x+\epsilon.y:x,y\in\mathbb{Z}_{2k},\epsilon^{2}=0\}
\]
in order to attach to each \emph{cantus firmus} $x$ the interval $y$ that separates
it from its \emph{discantus}\footnote{The discantus can be understood in the \emph{sweeping}
($x+y$) or the \emph{hanging} ($x-y$) orientations, but we will only use the sweeping
orientation from this point on.}. Thus for a strong dichotomy $\Delta=(X/Y)$ we have the consonant intervals
\[
 X[\epsilon]:=\{c+\epsilon.x:c\in\mathbb{Z}_{2k},x\in X\}
\]
and the dissonant intervals $Y[\epsilon]=\mathbb{Z}_{2k}\setminus X[\epsilon]$. Considering
the group
\[
\overrightarrow{GL}(\mathbb{Z}_{2k}[\epsilon]):=\{T^{a+\epsilon.b}.(v+\epsilon.w):a,b,w\in\mathbb{Z}_{2k},v\in\mathbb{Z}_{2k}^{\times}\},
\]
there is a canonical autocomplementary symmetry $p_{\Delta}^{c}\in\overrightarrow{GL}(\mathbb{Z}_{2k})$ such that
\[
 p_{\Delta}^{c}(X[\epsilon]) = Y[\epsilon]\quad\text{and}
\]
and leaves the \emph{tangent space} $c+\epsilon.\mathbb{Z}_{2k}$ invariant.

With this preamble it is possible to state a classical paradox for first-species counterpoint theory:
all the intervals $c+\epsilon.k$ used in a first-species counterpoint composition or improvisation are consonances.
Hence, how can any tension between the voices arise, if at all? Mazzola's solution is inspired in the fact \cite[p. 33-35]{kjS74}
that it is not that the point $c$ which is to be confronted against $c+k$, but it is the consonant point $\xi=c_{1}+\epsilon.k_{1}$
who will face a successor $\eta=c_{2}+\epsilon.k_{1}$. The idea is to \emph{deform} the dichotomy $(X[\epsilon]/Y[\epsilon])$ into
$(gX[\epsilon],gY[\epsilon])$ through a symmetry
$g\in\overrightarrow{GL}(\mathbb{Z}_{2k}[\epsilon])$, such that
\begin{enumerate}
\item the interval $\xi$ becomes a deformed dissonance, i.e., $\xi\in gY[\epsilon]$,
\item the symmetry $g$ is an autocomplementary function of
\[
 (gX[\epsilon],gY[\epsilon])
\]
which means that $p(gX[\epsilon])=gY[\epsilon]$,
\end{enumerate}
and thus we can transit from $\xi$ to a consonance $\eta$ which is also a deformed consonance, i.e., $\eta\in gX[\epsilon]\cap X[\epsilon]$. Since we wish to have
the maximum amount of choices, we request also that
\begin{enumerate}[resume]
\item the set $gX[\epsilon]\cap X[\epsilon]$ is of maximum cardinality among the symmetries that satisfy conditions 1 and 2.
\end{enumerate}
The elements of this latter set are the \emph{admitted successors}.

\section{Dichotomies of $2$-intervals}

For the purposes of the second-species counterpoint, we need now an
algebraic structure such that two intervals can be attached to a base
tone. In the spirit of the model presented in the previous section, we take all the polynomials of the
form
\[
 c+\epsilon_{1}.x+\epsilon_{2}.y\in \frac{\mathbb{Z}_{2k}[\mathcal{X},\mathcal{Y}]}{\langle \mathcal{X}^{2},\mathcal{Y}^{2},\mathcal{XY}\rangle}
 = \mathbb{Z}_{2k}[\epsilon_{1},\epsilon_{2}]
\]
where $\epsilon_{1}\equiv \mathcal{X} \bmod\langle \mathcal{X}^{2},\mathcal{Y}^{2},\mathcal{XY}\rangle$,
$\epsilon_{2}\equiv \mathcal{Y} \bmod \langle \mathcal{X}^{2},\mathcal{Y}^{2},\mathcal{XY}\rangle$, $c$ is the cantus firmus and $x,y$ are the intervals ($x$
is for the downbeat and $y$ is for the upbeat). An element $\xi\in\mathbb{Z}_{2k}[\epsilon_{1},\epsilon_{2}]$ is called a \emph{$2$-interval}. If
$\Delta=(X/Y)$ is a strong dichotomy with
polarity $p=T^{u}\circ v$, then
\[
 X[\epsilon_{1},\epsilon_{2}]:=\mathbb{Z}_{2k}+\epsilon_{1}.X+\epsilon_{2}.\mathbb{Z}_{2k}
\]
is an dichotomy in $\mathbb{Z}_{2k}[\epsilon_{1},\epsilon_{2}]$. We choose this dichotomy because the rules of counterpoint
demand that the interval that comes on the downbeat to be a consonance. A polarity for this
dichotomy, which is analogous to the one for the first-species case, is
\[
 p^{c}=T^{c(1-v)+\epsilon_{1}.u+\epsilon_{2}.u}\circ v
\]
because
\begin{align*}
 p^{c}X[\epsilon_{1},\epsilon_{2}] &= T^{c(1-v)}\circ v.\mathbb{Z}_{2k}+\epsilon_{1}.pX + \epsilon_{2}.p\mathbb{Z}_{2k}\\
 &= \mathbb{Z}_{2k}+\epsilon_{1}.Y+\epsilon_{2}.\mathbb{Z}_{2k}\\
 &= Y[\epsilon_{1},\epsilon_{2}]
\end{align*}
and it is such that
\[
 p^{c}(c+\epsilon_{1}.\mathbb{Z}_{2k}+\epsilon_{2}.\mathbb{Z}_{2k}) = c+\epsilon_{1}.\mathbb{Z}_{2k}+\epsilon_{2}.\mathbb{Z}_{2k},
\]
which means $p^{c}$ fixes the tangent space to cantus firmus $c$ as well.

We also check the following formula for future use:
\begin{align*}
 p^{c_{1}+c_{2}} &= T^{(c_{1}+c_{2})(1-v)+\epsilon_{1}.u+\epsilon_{2}.u}\circ v\\
&= T^{c_{1}(1-v)+c_{2}(1-v)+\epsilon_{1}.u+\epsilon_{2}.u}\circ v\\
&= T^{c_{1}}\circ T^{-vc_{1}}\circ T^{c_{2}(1-v)+\epsilon_{1}.u+\epsilon_{2}.u}\circ v\\
&= T^{c_{1}}\circ T^{c_{2}(1-v)+\epsilon_{1}.u+\epsilon_{2}.u}\circ v\circ T^{-c_{1}}\\
&= T^{c_{1}}\circ p^{c_{2}} \circ T^{-c_{1}}.
\end{align*}

\section{Species Projections}
If we represent the polynomial $c+\epsilon_{1}.x+\epsilon_{2}.y$ as a
column vector, the the candidates to (non-invertible) \emph{species projections} are
\begin{align*}
 g:\mathbb{Z}_{2k}[\epsilon_{1},\epsilon_{2}]&\to \mathbb{Z}_{2k}[\epsilon_{1}]\\
 \begin{pmatrix}c\\x\\y\end{pmatrix}&\mapsto \begin{pmatrix}
s & 0 & 0 \\
sw_{1} & s & sw_{2}
\end{pmatrix}
\begin{pmatrix}c\\x\\y\end{pmatrix}
+
\begin{pmatrix}t_{1}\\t_{2}\end{pmatrix}\\
&\phantom{=}= [sc+t_{1}]+\epsilon_{1}.[s(w_{1}c+x+w_{2}y)+t_{2}]
\end{align*}
for we want to keep it as simple as possible and that the second part of the interval to influence the
first part of the successor, but not the second one. We do
not require the transformation to be bijective for we want it to
be able to swap from second-species to first-species if
necessary\footnote{For the converse swap the
standard rules of counterpoint suffice: we can arbitrarily define the
third component of the $2$-interval. This is coherent with the
local application of counterpoint rules in Fux's theory, and also with the
particular idea of projection that stems from the fact that, in order to analyze
a fragment, we ``disregard'' notes on the upbeat \cite[pp. 41-43]{aM65}.}.

Let $X[\epsilon_{1},\epsilon_{2}.y] :=
\mathbb{Z}_{12}+\epsilon_{1}.X+\epsilon_{2}.y$. We might define a
species projection of a $2$-interval $\xi = c+\epsilon_{1}.x+
\epsilon_{2}.y$ as one such that
\begin{enumerate*}
\item the condition $c+\epsilon_{1}.x \notin gX[\epsilon_{1},\epsilon_{2}.y]$ holds,
\item the square
\begin{equation}\label{E:Comm}
\begin{CD}
\mathbb{Z}_{2k}[\epsilon_{1},\epsilon_{2}] @>g>> \mathbb{Z}_{2k}[\epsilon_{1}]\\
@Vp^{c}VV @VVp^{c}_{\Delta}V\\
\mathbb{Z}_{2k}[\epsilon_{1},\epsilon_{2}] @>>g> \mathbb{Z}_{2k}[\epsilon_{1}]
\end{CD}
\end{equation}
commutes, where
\[
 p_{\Delta}^{c}:=T^{c(1-v)+\epsilon_{1}.u}\circ v
\]
is the \emph{canonical} polarity of $(X[\epsilon_{1}]/Y[\epsilon_{1}])$, and
\item the cardinality of $gX[\epsilon_{1},\epsilon_{2}.y]\cap X[\epsilon_{1}]$
is maximal among the projections with the previous properties.
\end{enumerate*}

The reason for the second requirement is that when it is fulfilled then
\[
 p^{c}_{\Delta}(gX[\epsilon_{1},\epsilon_{2}]) = g(p^{c}X[\epsilon_{1},\epsilon_{2}])
=gY[\epsilon_{1},\epsilon_{2}],
\]
thus $p^{c}_{\Delta}$ is an autocomplementary function of $gX[\epsilon_{1},\epsilon_{2}]$.

\section{Algorithm for the Calculation of Projections}

As with the first-species case, if for a projection of the form
\[
g = T^{\epsilon_{1}.t_{2}}\circ \begin{pmatrix} s & 0 & 0 \\ sw_{1} & s & sw_{2}\end{pmatrix}
\]
we define
\[
 g^{(t_{1})} = g\circ T^{\epsilon_{1}{.}s^{-1}w_{1}t_{1}+\epsilon_{2}{.}t_{1}}
\]
then the relation
\[
 T^{t_{1}}\circ g = g^{(-t_{1})}\circ T^{s^{-1}t_{1}+\epsilon_{2}.t_{1}},
\]
holds, and hence contrapuntal projections can be calculated with cantus firmus
$0$ and successors can be suitably adjusted \cite[Theorem 2.2]{AJM15}. Therefore, we can set
$t_{1}=0$ and work with intervals of the form $\xi = \epsilon_{1}.y+\epsilon_{2}.z$.
For \eqref{E:Comm} to commute, it is necessary and sufficient that
\begin{equation}\label{E:Comm2}
 t_{2}+su(1+w_{2}) = u+vt_{2}.
\end{equation}

For $\epsilon_{1}.y \notin gX[\epsilon_{1},\epsilon_{2}.z]$ we need
\[
 y = sp(\ell)+t_{2}+sw_{2}z
\]
for some $\ell\in X$. Hence, for some $\ell\in X$ we have
\begin{equation}\label{E:CondS}
 t_{2} = y-s(p(\ell)+w_{2}z).
\end{equation}

\begin{remark}
Letting $w_{2}=0$ in \eqref{E:Comm2} and \eqref{E:CondS}, they reduce to the
first-species case. Thus, taking $s=v$ and $\ell = y$ both are satisfied and
hence we conclude that there exists at least one second-species counterpoint projection.
\end{remark}

We only need to work with the following set
\begin{align*}
gX[\epsilon_{1},\epsilon_{2}{.}z] &= \bigcup_{x\in\mathbb{Z}_{k}} g(x+\epsilon_{1}.X+\epsilon_{2}{.}z)\\
&= \bigcup_{x\in\mathbb{Z}_{2k}} (sx+\epsilon_{1}.(sw_{1}x+sw_{2}z+t_{2}+sX))\\
&= \bigcup_{r\in \mathbb{Z}_{2k}} (r+\epsilon_{1}.(w_{1}r+sX+w_{2}sz+t_{2}))\\
&= \bigcup_{r\in\mathbb{Z}_{2k}} (r+\epsilon_{1}.T^{w_{1}r+w_{2}sz+t_{2}}\circ s X)
\end{align*}
to calculate the following cardinality
\[
|gX[\epsilon_{1},\epsilon_{2}.z]\cap X[\epsilon_{1},\epsilon_{2}.z]| =\sum_{r\in\mathbb{Z}_{2k}}
|T^{w_{1}r+w_{2}sz+t_{2}}\circ s X \cap X|.
\]

When \eqref{E:CondS} holds, this reduces to
\begin{equation}\label{E:InterSE}
|gX[\epsilon_{1},\epsilon_{2}.z]\cap X[\epsilon_{1},\epsilon_{2}.z]| =\sum_{r\in\mathbb{Z}_{2k}}
|T^{w_{1}r+y-sp(\ell)}\circ s X \cap X|.
\end{equation}

From now on we only need to adapt \emph{mutatis mutandis}
Hichert's algorithm \cite[Algorithm 2.1]{AJM15} to search projections that maximize the intersection.

We must remark that \eqref{E:Comm2} and \eqref{E:CondS} are
perturbations of the conditions to find the counterpoint symmetries for
the first-species case. These, together with \eqref{E:InterSE},
show that the conditions for deducing a counterpoint theorem
\cite[Theorem 2.3]{AJM15} hold again, which yields the following result.

\begin{theorem}
Given a marked strong dichotomy $(X/Y)$ in $\mathbb{Z}_{2k}$, the $2$-interval $\xi\in X[\epsilon_{1},\epsilon_{2}]$
has at least $k^{2}$ and at most $2k^{2}-k$ admitted successors.
\end{theorem}

\begin{algo}
Here $\chi(x,y)$ is the function that returns the cardinality $T^{x}. y X \cap X$.
\begin{algorithmic}[1]
\REQUIRE A strong dichotomy $\Delta =(X/Y)$ and its polarity $T^{u}.v$.
\ENSURE The set of counterpoint projections $\Sigma_{y,z}\subseteq H$ for each
$\epsilon{.}y+\epsilon{.}z\in X[\epsilon_{1},\epsilon_{2}]$.
\FORALL{$y,z\in X$}
 \STATE $M\gets 0, \Sigma_{y,z} \gets \emptyset$.
 \FORALL{$s\in GL(\mathbb{Z}_{2k})$}
  \FORALL{$\ell\in X$}
   \FORALL{$w_{1},w_{2}\in\mathbb{Z}_{2k}$}
    \STATE $t_{2}\gets y-s((v\ell+u)+w_{2}z)$.
    \IF{$t_{2}+su(1+w_{2}) = u+vt_{2}$}
     \IF{$w_{1}=0$}
      \STATE $S\gets 2k\chi(t_{2},s)$.
     \ELSIF{$w_{1}\in GL(\mathbb{Z}_{2k})$}
      \STATE $S\gets k^{2}$
     \ELSE
      \STATE $\rho \gets \gcd(w_{1},2k)$
      \STATE $S\gets \rho\sum_{j=0}^{\frac{2k}{\rho}-1}\chi(j\rho+t_{2}+w_{2}z,s)$.
     \ENDIF
     \IF{$S>M$}
      \STATE $\Sigma_{y,z} \gets \left\{T^{\epsilon_{2}{.}t_{2}}\circ \begin{pmatrix}
s & 0 & 0 \\
sw_{1} & s & sw_{2}
\end{pmatrix}\right\}$.
      \STATE $S\gets M$.
     \ELSIF{$S=M$}
      \STATE $\Sigma_{y,z} \gets \Sigma_{y,z} \cup \left\{T^{\epsilon{.}t_{2}}\circ \begin{pmatrix}
s & 0 & 0 \\
sw_{1} & s & sw_{2}
\end{pmatrix}\right\}$.
     \ENDIF
    \ENDIF
   \ENDFOR
  \ENDFOR
 \ENDFOR
 \RETURN $\Sigma_{y,z}$.
\ENDFOR
\end{algorithmic}
\end{algo}

\begin{example}
The first (valid\footnote{The first example is the student's attempt to write a second-species discantus by himself, but he makes two mistakes
near the end of the exercise, namely the steps from the sequence $7+\epsilon_{1}.2+\epsilon_{2}.11$, $5+\epsilon_{1}.0+\epsilon_{2}.9$,
$4+\epsilon_{1}.11+\epsilon_{2}.1$. They are also forbidden steps in the projection model!}) example of second-species counterpoint in
the \emph{Gradus ad Parnassum} \cite[p. 45]{aM65} is
\begin{multline*}
 \xi_{1} = 2 + \epsilon_{1}.7 + \epsilon_{2}{.}0,\, \xi_{2} = 5 +
\epsilon_{1}{.}4 + \epsilon_{2}{.}6,\, \xi_{3} = 4 + \epsilon_{1}{.}8 +
\epsilon_{2}{.}3,\\ \xi_{4} = 2 + \epsilon_{1}{.}7 + \epsilon_{2}{.}0,\,
\xi_{5} = 7 + \epsilon_{1}{.}4 + \epsilon_{2}{.}5,\, \xi_{6} = 5 +
\epsilon_{1}{.}9 + \epsilon_{2}{.}4,\\ \xi_{7} = 9 + \epsilon_{1}{.}3 +
\epsilon_{2}{.}5,\, \xi_{8} = 7 + \epsilon_{1}{.}9 + \epsilon_{2}{.}4,\,
\xi_{9} = 5 + \epsilon_{1}{.}9 + \epsilon_{2}{.}4,\\ \xi_{10} = 4 +
\epsilon_{1}{.}7 + \epsilon_{2}{.}9,\,\xi_{11} = 2 + \epsilon_{1}{.}0
\end{multline*}

\begin{center}
 \includegraphics[width=\textwidth]{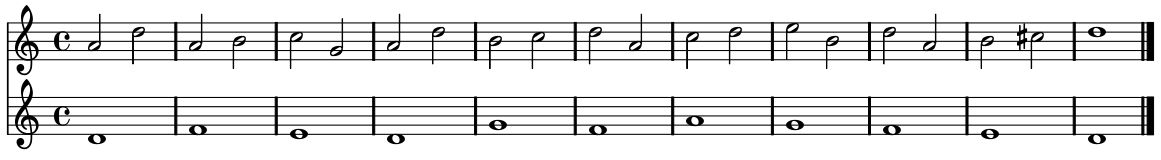}
\end{center}

Some counterpoint projections for the successors are
\begin{multline*}
 g_{1} = \begin{pmatrix}
7 & 0 & 0\\
0 & 7 & 0
\end{pmatrix},
 g_{2} = T^{\epsilon_{1}{.}6}\circ\begin{pmatrix}
1 & 0 & 0 \\
6 & 1 & 0
\end{pmatrix},
 g_{3} = T^{\epsilon_{1}{.}6}\circ\begin{pmatrix}
7 & 0 & 0\\
6 & 7 & 0
\end{pmatrix}\\
g_{4} = g_{1},
g_{5} = g_{2},
g_{6} = T^{\epsilon_{1}{.}8}\circ\begin{pmatrix}
5 & 0 & 0 \\
8 & 5 & 0
\end{pmatrix},\\
g_{7} = \begin{pmatrix}
11 & 0 & 0\\
0 & 11 & 8
\end{pmatrix},
g_{8} = g_{6},
g_{9} = g_{6},
g_{10} = g_{1}.
\end{multline*}

Let us examine in little bit more of detail the first transition. Note that $\eta=11+\epsilon_{1}.4+\epsilon_{1}.11$ is
a consonance, and that
\[
 g_{1}(\eta)=\begin{pmatrix}
 7 & 0 & 0 \\
 0 & 7 & 0
 \end{pmatrix}
 \begin{pmatrix}
 11\\
 4\\
 11
 \end{pmatrix}
 =\begin{pmatrix}
 5\\
 4
 \end{pmatrix},
\]
which justifies the fact that the $2$-interval $5+\epsilon_{1}.4+\epsilon_{2}.6$ is an admitted successor.
\end{example}


\section{Comparison with Fux's Approach}

Fux states the following in relation to second-species counterpoint (emphasis is our own) \cite[p. 41]{aM65}:
\begin{quote}
 The second species results when two half notes are set against a whole note. The first of
 them comes on the downbeat and must always be consonant; the second comes on the upbeat
 and \emph{it may be dissonant if it moves from the preceding note and to the following
 note stepwise}. However, \emph{if it moves by a skip, it must be consonant}.
\end{quote}

We made a program that compares the performance
of a first-species model that takes into account Fux's restrictions against the projection model. More explicitly,
taking a second-species step
\[
 (0+\epsilon_{1}.k_{1}+\epsilon_{2}.t_{1},c_{2}+\epsilon_{1}.k_{2})
\]
such that we can proceed (in first-species) from $0+\epsilon_{1}.k_{1}$ to $c_{2}+\epsilon_{1}.k_{2}$, we verify the following cases:
\begin{enumerate}
\item the upbeat interval $t_{1}$ of the first $2$-interval is allowed to be dissonant only when it connects a valid progression of consonances stepwise, i.e., $0+t_{1}$ is between $0+k_{1}$ and $c_{2}+k_{2}$ and it is separated at most $2$ semitones from them and
\item if $t_{1}$ is consonant, we duplicate
the cantus firmus and check if $(0+\epsilon.k_{1},0+\epsilon.t_{1})$ and $(0+\epsilon.t_{1},c_{2}+\epsilon.k_{2})$ are valid
first-species steps.
\end{enumerate}

The results appear in Table \ref{T:Comparacion} for cases 1 and 2. We must stress that the projection model was not restricted
in case 1 to stepwise dissonances but it allowed any dissonance in the upbeat.

\begin{table}[ht]
\begin{tabular}{|l|l|l|}
\hline
Number of steps                     & Case 1 & Case 2 \\ \hline
Total                               & 1994   & 2592   \\ \hline
Valid only for Fux model            & 9      & 178    \\ \hline
Valid only for the projection model & 1447   & 860    \\ \hline
Valid in both models                & 301    & 1464   \\ \hline
\end{tabular}
\caption{Data for comparison of Fux's model with restrictions for second species against the projection model.}
\label{T:Comparacion}
\end{table}

We note that the number of cases the projection model cannot explain and only Fux can is relatively small: they amount to $2.9$\%
and $17.1$\% for cases 1 and 2, respectively. Thus we can conclude that the vast majority of what is forbidden in the projection
model is also forbidden in Fux's model, or that we have successfully extended Fux's handling of dissonance and consonance for second species.
Even if this could be ascribed to the fact that the projection model admits $87.663$\% and $89.660$\% of the total of transitions in cases
1 and 2, respectively, it should be kept in mind that the one-species model admits $89.671$\% of the possible steps between consonant
intervals \cite[p. 48]{aN10}.

\bibliographystyle{amsplain}
\bibliography{2ndspecies}

\end{document}